\documentclass[12pt]{amsart}		
\usepackage[margin=1in]{geometry}	
\usepackage{pinlabel}
\usepackage{graphicx}
\usepackage{verbatim, color, footmisc}
\usepackage{hyperref}
\usepackage{amsthm}
\usepackage[all]{xy}
\usepackage{mathdots}
 \def\co{\colon\thinspace}
 \def\CC{{C\nolinebreak[4]\hspace{-.05em}\raisebox{.4ex}{\tiny\bf ++}}} 

\author{Thomas Peters }
\address{Department of Mathematics, Columbia University, MC 4406\\ 2990 Broadway, New York,  NY 10027}
\email{tpeters@math.columbia.edu}

\title{\textbf{A concordance invariant from the Floer homology of $\pm1$-surgeries}}
\begin{document}
\newtheorem{theorem}{Theorem}[section]
\newtheorem{lemma}[theorem]{Lemma}
\newtheorem{definition}[theorem]{Definition}
\newtheorem{proposition}[theorem]{Proposition}
\newtheorem{claim}[theorem]{Claim}
\newtheorem{rem}[theorem]{Remark}
\newtheorem{corollary}[theorem]{Corollary}

\begin{abstract}
We discuss a concordance invariant constructed from Heegaard Floer homology ``correction terms'' and $\pm1$ surgeries on knots. 
\end{abstract}

\maketitle


\section{Introduction}
\label{section:dIntro}
Given a closed oriented three-manifold with torsion $\mathrm{Spin}^c$ structure, the associated Heegaard Floer homology groups come with absolute $\mathbb{Q}$--gradings; see Ozsv\'ath--Szab\'o \cite{OSzFour}.  This allows one to define numerical invariants of $\mathrm{Spin}^c$ three-manifolds, the so-called ``correction terms'' or ``$d$--invariants''.  Specifically, suppose $(Y,\mathfrak{s})$ is a $\mathrm{Spin}^c$ rational homology three-sphere.  Then Ozsv\'ath and Szab\'o define $d(Y,\mathfrak{s})$ (the correction term) to be the minimal degree of any non-torsion class in $HF^+(Y,\mathfrak{s})$ coming from $HF^\infty(Y,\mathfrak{s})$\footnote{There are correction terms for three-manifolds with positive first Betti number, but we do not discuss them at the moment.}.  This invariant is analogous to the monopole Floer homology $h$--invariant introduced by Fr{\o}yshov \cite{Froy1}.  If $Y$ only has a single $\mathrm{Spin}^c$ structure $\mathfrak{s}_0$ (ie if $Y$ is an integer homology sphere), then we denote $d(Y,\mathfrak{s}_0)$ by just $d(Y)$.   The $d$--invariants satisfy some useful properties, according to the following theorem of Ozsv\'ath and Szab\'o:
\begin{theorem}[Ozsv\'ath--Szab\'o \cite{OSzAbsGr}]
\label{theorem:dprops}
Let $Y$ be an oriented rational homology three-sphere.  Its correction terms satisfy:
\begin{enumerate}
\item \label{item:conj} Conjugation invariance \[d(Y,\mathfrak{s}) = d(Y,\bar{\mathfrak{s}}).\]
\item \label{item:negDef} If $Y$ is an integral homology three-sphere and is the oriented boundary of a negative-definite four-manifold $W$ then $d(Y)\geq0$.
\end{enumerate}
\end{theorem}
In fact, item \ref{item:negDef} follows from a more general statement, Proposition \ref{prop:ineq}, and the following theorem of Elkies.
\begin{theorem}[Elkies \cite{Elk}]
\label{theorem:elkies}  Let $Q\co V\otimes V\to\mathbb{Z}$ be a negative-definite unimodular bilinear form over $\mathbb{Z}$.  Denote by $\Xi(Q)$ the set of \emph{characteristic vectors for Q}, ie the set of vectors $\xi\in V$ satisfying $$Q(\xi,v)\equiv Q(v,v)\mod2$$ for all $v\in V$.  Then, $$0\leq\max_{\xi\in \Xi(Q)} Q(\xi,\xi) + \mathrm{dim}(V),$$ with equality if and only if the bilinear form $Q$ is diagonalizable over $\mathbb{Z}$.
\end{theorem}
Also, $Y$ can be a disjoint union of rational homology three-spheres, in which case Theorem \ref{theorem:dprops} (together with Theorem \ref{theorem:elkies}) implies:
\begin{corollary}[Ozsv\'ath--Szab\'o \cite{OSzAbsGr}]
\label{corollary:addendum}
Let $Y_1$ and $Y_2$ be oriented rational homology three-spheres.  Then
\begin{enumerate}
\item \label{item:reverse} Let $-Y_1$ denote the manifold $Y_1$ with opposite orientation, then $$d(Y_1,\mathfrak{s}) = -d(-Y_1,\mathfrak{s})$$
\item \label{item:cobord} If $(Y_1,\mathfrak{s}_1)$ is $\mathrm{Spin}^c$ rational homology cobordant to $(Y_2,\mathfrak{s}_2)$, then \[d(Y_1,\mathfrak{s}_1) = d(Y_2,\mathfrak{s}_2).\]
\item \label{item:negDefCob} If $Y_1$ and $Y_2$ are integral homology three-spheres and $W$ is a negative-definite $\mathrm{Spin}^c$ cobordism from $(Y_1,\mathfrak{s}_1)$ to $(Y_2,\mathfrak{s}_2)$, then \[d(Y_2,\mathfrak{s}_2) \geq d(Y_1,\mathfrak{s}_1).\]  
\item \label{item:QB4} If $(Y_1,\mathfrak{s}_1)$ bounds a rational homology four-ball, then $d(Y,\mathfrak{s}) = 0$.
\end{enumerate}
\end{corollary}

Heegaard Floer homology $d$--invariants have been used to give restrictions on intersection forms of four-manifolds which can bound a given three-manifold (for instance, Ozsv\'ath and Szab\'o reproved Donaldson's diagonalization theorem using correction terms).  They have also been used to define concordance invariants of knots\footnote{Heegaard Floer theory has led to other concordance invariants, most notably the $\tau$ invariant (see Ozsv\'ath--Szab\'o \cite{OSTau} and Rasmussen \cite{Ra}).}.  For instance, Manolescu and Owens \cite{MO} used the $d$--invariants of the branched double cover of a knot to produce concordance invariants (see also Grigsby, Ruberman, and Strle \cite{GRS}, Jabuka \cite{Jab2}, and Jabuka and Naik \cite{Jab}).  In this paper, given a knot $K\subset S^3$ in the three-sphere, we show that $d(S^3_{+1}(K))$ is a concordance invariant of $K$ and examine some of its properties.
 Occasionally we denote this invariant by $dS^3_1$.  Note that one could also study $d(S^3_{-1}(K))$, but these invariants are determined by the $dS^3_1$ since $d(S^3_{-1}(K)) = d(-S^3_{1}(mK)) = -d(S^3_{1}(mK))$ where $mK$ denotes the mirror of $K$.  We also establish a ``skein inequality'' reminiscent of a property of the knot signature.  Specifically,
\begin{theorem} Given a diagram for a knot with distinguished crossing $c$, let $D_+$ and $D_-$ be the result of switching $c$ to positive and negative crossings, respectively, as in Figure \ref{figure:skein}.  Then
\label{theorem:dIneq} \[ d(S^3_{1}(D_-);\mathbb{F}) - 2 \leq d(S^3_{1}(D_+);\mathbb{F}) \leq d(S^3_{1}(D_-);\mathbb{F}) \] for any field $\mathbb{F}$.  Here, $d(Y;\mathbb{F})$ denotes the correction term of $Y$ computed from Floer homology with coefficients in $\mathbb{F}$.
\end{theorem}Indeed, we expect that the restriction that the coefficients are taken in a field could be relaxed to include $\mathbb{Z}$--coefficients, but our proof only holds for field coefficients.  The invariants $d(S^3_{\pm1}(K))$ also give rise to four-ball genus bounds.  Specifically, we have the following:
\begin{theorem} 
\label{theorem:genus} Let $K$ be a knot in the three-sphere.  Then
\[
	0 \leq -d(S^3_1(K);\mathbb{Z}_2) \leq 2 g_4(K)
\] where $g_4(K)$ denotes the smooth four-ball genus of $K$.
\end{theorem}

Again, we expect that this should hold for Floer homology with any coefficients, but our proof is special to $\mathbb{Z}_2$--coefficients.  Theorem \ref{theorem:genus} should be compared to the following theorem of Fr{\o}yshov:
\begin{theorem}[Fr{\o}yshov \cite{Froy2}]
\label{theorem:froy-ineq} Let $Y$ be an oriented homology three-sphere and $\gamma$ a knot in $Y$ of ``slice genus'' $\widetilde{g}$.  If $Y_{\gamma,-1}$ is the result of $-1$--surgery on $\gamma$ then $$0\leq h(Y_{\gamma,-1}) - h(Y) \leq \lceil \widetilde{g}/2\rceil.$$
\end{theorem}
Here $h(Y)$ is Fr{\o}yshov's instanton Floer homology $h$--invariant and the ``slice genus'' is defined to be the smallest non-negative integer $\widetilde{g}$ for which there exists a smooth rational homology cobordism $W$ from $Y$ to some rational homology sphere $Y'$ and a genus $\widetilde{g}$ surface $\Sigma\subset W$ such that $\partial W = \gamma$.  It is not clear to the author whether this definition agrees with the usual one for $Y = S^3$.  In light of the conjectural relationship $h(Y) = d(Y)/2$ and Theorem \ref{theorem:froy-ineq}, we suspect that the inequality in Theorem \ref{theorem:genus} is in general weaker than the $h$--invariant inequality.

Finally, using the theory of Ozsv\'ath--Szab\'o \cite{OSzKnot} and Rasmussen \cite{Ra}, we observe how one can algorithmically compute $d(S^3_{\pm1}(K))$ if one knows the filtered chain homotopy type of the knot complex $CFK^\infty(K)$.  A computer implementation of this algorithm is discussed.  

\subsection{Further questions}
What is the relationship between the correction terms of $\pm1$--surgeries on a knot and the Ozsv\'ath--Szab\'o, Rasmussen $\tau$ invariant?  From the discussion in Section \ref{section:computations}, it seems likely that $|d(S^3_{\pm1}(K)|\leq 2|\tau(K)|$, but as of the time of this writing a proof remains elusive.  Of course if this were the case, then the genus bound, Theorem \ref{theorem:genus}, would follow immediately from the inequality $|\tau(K)| \leq g_4(K)$ (see Ozsv\'ath and Szab\'o \cite{OSTau} for a discussion).

\subsection{Organization}

This paper is organized as follows.  In Section \ref{section:theInvariant}, we discuss basic properties of $d(S^3_1(K))$, including its invariance under concordance.  In Section \ref{section:skein} we give a proof of the skein inequality, Theorem \ref{theorem:dIneq}.  In Section \ref{section:genus} we prove Theorem \ref{theorem:genus}.  Finally in Section \ref{section:computations} we discuss an algorithm to compute $d(S^3_1(K))$ given the knot complex $CFK^\infty(K)$ as well as a computer implementation of this algorithm.

\subsection{Acknowledgement}

The author would like to thank his PhD supervisor, Peter Ozsv\'ath for suggesting the problem as well as invaluable guidance over the years.  He would also like to thank Maciej Borodzik, Kim Fr{\o}yshov, Matt Hedden, Adam Levine, and Danny Ruberman for helpful conversations.

\section{The invariant}
\label{section:theInvariant}

\begin{figure}[ht!]
\labellist
\small\hair 2pt
\pinlabel $A$ at 56 82
\pinlabel $K_1$ at -8 37
\pinlabel $W$ at 60 -8
\pinlabel $h^2$ at 126 11
\pinlabel $K_2$ at 150 67
\endlabellist
\begin{center}
\includegraphics{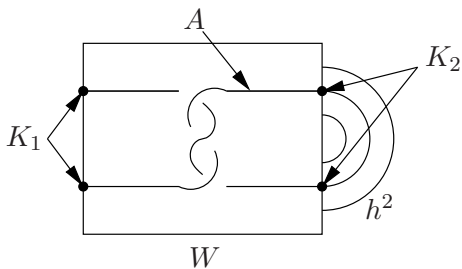}
\caption{\label{figure:knotConcord}The cobordism $W$.    }
\end{center}
\end{figure}

\begin{proposition} $d(S^3_{1}(K))$ is a concordance invariant.
\end{proposition}
\begin{proof}
It is simple to see that $d(S^3_{1}(K)) = 0$ if $K$ is smoothly slice: $S^3_{1}(K)$ bounds the four-manifold obtained by attaching a $+1$--framed two-handle along $K$ to the four-ball.  This four-manifold has second homology generated by a sphere of square $+1$.  By blowing this down, we see that $S^3_{1}(K)$ bounds a rational homology four-ball.  By item \ref{item:QB4} of Corollary \ref{corollary:addendum}, it follows that $d(S^3_{1}(K)) = 0$.  It is just slightly more work to see that if $K_1$ and $K_2$ are smoothly concordant, then $d(S^3_{1}(K_1)) = d(S^3_{1}(K_2))$:  the concordance gives us a smoothly embedded annulus $A \subset S^3\times I$ (here $I=[0,1]$) with $\partial A = K_1 \cup K_2$ and $K_1\subset S^3\times \{0\}$, $K_2\subset S^3\times\{1\}$.  Attach a two-handle to $S^3\times\{1\}$ with framing $+1$ along $K_2$ to give a four-manifold $W$ (see Figure \ref{figure:knotConcord}).  Consider a small regular neighborhood of the core disk of this two-handle union a regular neighborhood of the annulus $A$.  This gives cobordisms $W_0:S^3\to S^3_{1}(K_1)$ and $W_1:S^3_{1}(K_1)\to S^3_{1}(K_2)$
such that $W = W_0\cup W_1$.  Notice that $W_0$ is just $S^3\times [0,\varepsilon] \cup h^2$, a $+1$--framed two-handle attached along $K_1$ to a thickened $S^3$.  It follows that $b_2(W_0) = b_2(W) = 1$ and $b_2(W_1) = 0$ (this last fact can be seen from the Mayer--Vietoris sequence applied to the decomposition $W = W_0\cup W_1$: $0 = H_2(W_0\cap W_1;\mathbb{Z}) \to H_2(W_0;\mathbb{Z})\oplus H_2(W_1;\mathbb{Z}) \cong \mathbb{Z}\oplus H_2(W_1;\mathbb{Z}) \to H_2(W;\mathbb{Z}) \cong \mathbb{Z} \to H_1(W_0\cap W_1;\mathbb{Z}) = 0$).  Applying item \ref{item:cobord} of Corollary \ref{corollary:addendum}  to $W_1$ shows that $d(S^3_{1}(K_1)) = d(S^3_{1}(K_2))$.
\end{proof}

In fact, the basic topological fact that if knots $K$, $K'$ are concordant then $S^3_{\pm1}(K)$ is homology cobordant to $S^3_{\pm1}(K)$ used in the previous argument follows from a more general fact due to Gordon \cite{Gor1}: If two knots $K$ and $K'$ are concordant, then for any $r\in\mathbb{Q}$, we have a homology cobordism $S_r^3(K)\to S^3_r(K')$.  As pointed out by several people, this implies that for each rational $r\in\mathbb{Q}$, we get concordance invariants $d(S^3_r(K))$.  It is natural to ask about the independence of these invariants.

In general, calculating $d$--invariants is quite challenging.  However, in certain cases explicit formulae exist.  For instance, let $K$ be an alternating knot.  Then in \cite{OSalt}, Ozsv\'ath and Szab\'o prove that 
\begin{equation}
\label{eq:alt}
d(S^3_{+1}(K)) = 2\,\mathrm{min}\left(0,-\lceil \frac{-\sigma(K)}{4} \rceil\right)
\end{equation} where $\lceil x \rceil$ is the ceiling function and $\sigma(K)$ denotes the knot signature (see also Rasmussen \cite{RaTwo}).  This formula shows that the concordance invariants $d(S^3_{\pm1}(K))$ do not give group homomorphisms from the smooth concordance group to $\mathbb{Z}$: take the knot $RHT\,\#\,LHT$ where $RHT$ denotes the right-handed trefoil and $LHT$ denotes the left-handed trefoil.  This knot is slice and hence has vanishing $d S^3_{1}$ but $d(S^3_{+1}(RHT)) = -2$ and $d(S^3_{+1}(LHT)) = 0$.  Explicit formulae for $d$--invariants also exist in the case of certain plumbed three-manifolds; see Ozsv\'ath--Szab\'o \cite{OSPlumb}.  In another direction, since torus knots admit lens space surgeries, one may use Ozsv\'ath--Szab\'o \cite[Theorem 1.2]{OSLens} to calculate $dS^3_1$ for torus knots.

It may be worth noting that Equation \ref{eq:alt} does not hold for all knots.  For instance, the $(3,4)$--torus knot has signature $-6$ and $dS^3_1 = -2$.

The non-additivity of $dS^3_{+1}$ can be used to detect relations or establish linear independence in the smooth concordance group, $\mathcal{C}$.  For example, recall that $\sigma(LHT)=2$, $\tau(LHT)=-1$, and $s(LHT)=2$ (here, $s(K)$ denotes the Rasmussen $s$ concordance invariant of \cite{JakeS}).  It is also the case that $\sigma(T_{3,4}) = -6$, $\tau(T_{3,4}) = 3$, $s(T_{3,4}) = -6$ where here $T_{3,4}$ denotes the $(3,4)$--torus knot.  It follows that $f(LHT\,\#\,LHT\,\#\,LHT\,\#\,T_{3,4}) = 0$ for any $f$ among $s$, $\tau$, or $\sigma$.  However, this knot is not slice, since $d(S^3_{-1}(LHT\,\#\,LHT\,\#\,LHT\,\#\,T_{3,4})) = 2$, a fact which can be verified with our program \texttt{dCalc}.

Finally, note that $d(S^3_1(K))$ is always even.  This follows immediately from the long exact sequence
\begin{equation}
\label{equation:exact}
	\cdots \to HF^+(S^3) \to HF^+(S^3_0(K)) \to HF^+(S^3_1(K)) \to \cdots
\end{equation} of Ozsv\'ath and Szab\'o \cite{OSzAnn2}, and the fact that $HF^\infty(S^3_0(K))$ is standard (see, for instance, Section \ref{section:skein} for a discussion).

\section{Skein relations}
\label{section:skein}
\begin{figure}[ht!]
\labellist
\small\hair 2pt
\pinlabel $D_+$ at 19 -8
\pinlabel $D_-$ at 126 -8
\endlabellist
\begin{center}
\includegraphics{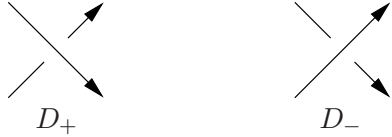}
\caption{  \label{figure:skein}  Positive and negative crossings, respectfully.}
\end{center}
\end{figure}

Recall the axiomatic characterization of the knot signature $\sigma$ found by Giller (see also Murasugi \cite{Murasugi}).  
\begin{theorem}[Giller \cite{Gil}]
\label{theorem:gil}
Suppose that $K$ is a knot (but not a link) and $D$ is a diagram for $K$.  Then $\sigma(K)$ can be determined from the following three axioms:
\begin{enumerate}
\item \label{item:M1} If $K$ is the unknot then $\sigma(K) = 0$.
\item \label{item:M2} If $D_+$ and $D_-$ are as in Figure \ref{figure:skein}, then \[ \sigma(D_-)-2 \leq \sigma(D_+) \leq \sigma(D_-)\] (recall that $\sigma$ is always even).
\item \label{item:M3} If $\Delta_K(t)$ is the Conway-normalized Alexander polynomial of $K$, then \[\mathrm{sign}(\Delta_K(-1)) = (-1)^{\sigma(K)/2}.\]   
\end{enumerate}
\end{theorem}

These axioms of course cannot hold for the invariant $d S^3_{1}$, but Theorem \ref{theorem:dIneq} does give us an analogue of Theorem \ref{theorem:gil}, item (2).

In light of the the axiomatic description of $\sigma$, it is an interesting question to calculate $dS^3_1/2$ modulo 2.  If one could achieve this, it might then be possible to give a completely algorithmic description of $dS^3_1$.

\begin{figure}
\labellist
\small\hair 2pt
\pinlabel $-1$ at 44 18
\pinlabel $\langle1\rangle$ at 45 46
\pinlabel $W_1$ at 18 -8
\pinlabel $+1$ at 151 18
\pinlabel $\langle1\rangle$ at 152 46
\pinlabel $W_0$ at 127 -8
\endlabellist
\begin{center}
\includegraphics{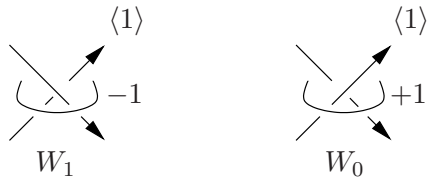}
\caption{ \label{figure:kirby1} A pair of relative handlebodies, representing the cobordisms $W_0$ and $W_1$. }
\end{center}
\end{figure}

\begin{figure}
\labellist
\small\hair 2pt
\pinlabel $-1$ at 86 32
\pinlabel $\langle1\rangle$ at 100 -8
\endlabellist
\begin{center}
\includegraphics{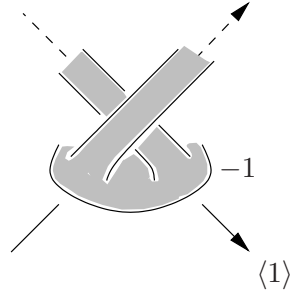}
\caption{ \label{figure:surf}The torus $T$ is represented by the shaded region, which is then capped off by the core of the $-1$--framed two-handle.   }
\end{center}
\end{figure}

We now return to the proof of Theorem \ref{theorem:dIneq}.
\begin{proof}[Proof of Theorem \ref{theorem:dIneq}]
{\bf Step 1: $d(S^3_1(D_-))\geq d(S^3_1(D_+))$:}

Given a knot $K$ with diagram $D(K)$ and a distinguished crossing, we have cobordisms $W_0\co S^3_1(D_-) \to S^3_{1}(D_+)$ and $W_1\co S^3_1(D_+) \to S^3_1(D_-)$ given by the Kirby diagrams in Figure \ref{figure:kirby1}.  We claim that $b_2(W_i)=1$ for $i=0,1$.  We argue this for $W_1$, the argument for $W_0$ being analogous.  $W_1$ fits into a four-manifold $W = X\cup W_1$ where $X$ is obtained by attaching a $+1$--framed two-handle along $K\subset S^3 = \partial B^4$ to the four-ball.  Clearly $b_2(W)=2$ and $b_2(X)=1.$  Consider the Mayer--Vietoris sequence applied to the decomposition $W = X\cup W_1$.  In this case, $X\cap W_1 \cong S^3_1(D_+)$, an integral homology three-sphere.  So we have $0 = H_2(X\cap W_1;\mathbb{Z}) \to H_2(X;\mathbb{Z})\oplus H_2(W_1;\mathbb{Z}) \cong \mathbb{Z}\oplus H_2(W_1;\mathbb{Z}) \to H_2(W) \cong \mathbb{Z}^2 \to H_1(X\cap W_1;\mathbb{Z}) = 0,$ showing that $H_2(W_1;\mathbb{Z}) \cong \mathbb{Z}$.  We may even find a torus $T$ in $W_1$ which generates $H_2(W_1;\mathbb{Z})\cong \mathbb{Z}$ as in Figure \ref{figure:surf}.  We claim that $[T]^2 = -1$ (likewise, for $W_0$ we have can find a torus of square $+1$ generating $H_2(W_1;\mathbb{Z})$).  $T$ of course sits inside the larger cobordism $W$.  Let $\{\alpha,\beta\}$ be an ordered basis of $H_2(W;\mathbb{Z})$ coming from the two two-handles (more specifically, $\alpha$ is the homology class of the core disk of the two-handle attached to $D_+$ capped off by a Seifert surface, and $\beta$ is the homology class of the core disk of the two-handle attached to the $-1$--framed knot in Figure \ref{figure:surf} capped off with a Seifert surface pushed slightly into the four-ball).  With respect to this basis, we see that the intersection form of $W$ is given by the matrix
\[
  Q_W = \begin{pmatrix} +1 & 0 \\
			0    & -1\\
	\end{pmatrix}.
\] By Figure \ref{figure:surf}, it is clear that $[T]\cdot\alpha = 0$ and $[T]\cdot\beta = -1$.  Therefore $[T] = \beta$ and $[T]^2 = -1$.  By item \ref{item:negDefCob} of Corollary \ref{corollary:addendum} it follows that 
\[
  d(S^3_1(D_-)) \geq d(S^3_1(D_+)).
\]

{\bf Step 2: $d(S^3_1(D_+)) \geq d(S^3_1(D_-)) - 2$:}\\
By a similar argument to the previous, we see that the cobordism $W_0\co S^3_1(D_-)\to S^3_1(D_+)$ has second homology generated by a torus $T$ of square $+1$.  Taking an internal connected sum  of $S^3_1(D_-)$ with a regular neighborhood of $T$, $\nu(T)$, we get cobordisms $V_0\co S^3_1(D_-)\to S^3_1(D_-)\,\#\,\partial\nu(T)$ and $V_1\co S^3_1(D_-)\,\#\,\partial\nu(T)\to S^3_1(D_+)$ such that $W_0 = V_0\cup V_1$.  Here $\partial\nu(T)$ denotes the boundary of a regular neighborhood of the surface $T$ in $W_0$.  This is of course a circle bundle over the two-torus with Euler number $+1$ (it is also a torus bundle over the circle with reducible monodromy).  It may be realized as $(0,0,1)$--surgery on the Borromean rings, which we denote by $M\{0,0,1\}$.  Clearly $b_2^+(V_1)=0$ for otherwise we would have a surface $F$ with positive square in $W_0$ which does not intersect the generating torus $[T]\in H_2(W_0;\mathbb{Z})$.  Similarly, we have that $b_2^-(V_1) = 0$.  Notice that $V_0$ deformation retracts onto the wedge $S^3_1(D_-)\,\vee\, T^2$ and hence has euler characteristic $\chi(V_0) = -1$.  Since $\chi(W_1) = \chi(V_0) + \chi(V_1) - \chi(V_0\cap V_1)$ and $\chi(V_0\cap V_1) = 0$ ($V_0\cap V_1$ is a three-manifold) we see that $\chi(V_1) = 2$.  Therefore the cobordism $V_1\co S^3_1(D_-)\,\#\,\partial\nu(T)\to S^3_1(D_+)$ has $\chi = 2$, $\sigma = 0$ (here $\sigma$ denotes the signature of the intersection form of $V_1$), and $c_1(\mathfrak{s})^2 = 0$ for all $\mathrm{Spin}^c$ structures $\mathfrak{s}\in \mathrm{Spin}^c\; V_1$.  By the formula for grading shifts in Heegaard Floer homology (see Ozsv\'ath and Szab\'o \cite{OSzFour}), it follows that the maps on Floer homology associated with this cobordism have grading shift
\[
  \frac{c_1(\mathfrak{s})^2 - 2\chi(V_1) - 3 \sigma(V_1)}{4} = -1.
\]

Before continuing with $d$--invariant calculations, we pause to recall some constructions in Heegaard Floer theory for manifolds with $b_1>0$.  In this case, there is a natural action of the exterior algebra $\Lambda^*H_1(Y;\mathbb{Z})/\mathrm{Tors}$ on all versions of Floer homology $HF^\circ$.  Under this action, elements of $H_1(Y;\mathbb{Z})$ drop relative gradings by one.  As an example, let $\mathbb{Z}_{(k)}$ denote the graded abelian group $\mathbb{Z}$ supported in grading $k$.  Under the graded isomorphism $\widehat{HF}(S^2\times S^1) \cong \mathbb{Z}_{(1/2)}\langle a\rangle \oplus \mathbb{Z}_{(-1/2)}\langle b\rangle$, the action of the circle factor $\gamma := [*\times S^1] \in H_1(S^2\times S^1;\mathbb{Z})$ is given by $\gamma\cdot a = b$ and $\gamma\cdot b = 0$.  

A $\mathrm{Spin}^c$ three-manifold $(Y,\mathfrak{s})$ with torsion $\mathrm{Spin}^c$ structure $\mathfrak{s}$ is said to have \emph{standard} $HF^\infty$ if there is a graded isomorphism of $\Lambda^*H_1(Y;\mathbb{Z})/\mathrm{Tors}\otimes_\mathbb{Z}\mathbb{Z}[U,U^{-1}]$--modules
\begin{equation}
\label{equation:std}
  HF^\infty(Y,\mathfrak{s}) \cong \Lambda\!^*H^1(Y;\mathbb{Z})\otimes_\mathbb{Z}\mathbb{Z}[U,U^{-1}]
\end{equation}
where the action of $H_1(Y;\mathbb{Z})/\mathrm{Tors}$ on the right hand side is given by contraction on\\ $\Lambda^*H^1(Y;\mathbb{Z})$.  Here $\Lambda^*H^1(Y;\mathbb{Z})$ is graded by the requirement that $\mathrm{gr}\left(\Lambda^{b_1Y}H^1(Y;\mathbb{Z})\right) = b_1(Y)/2$ and the fact that $H_1(Y;\mathbb{Z})$ drops gradings by 1.  For example, $\#^kS^2\times S^1$ has standard $HF^\infty$ for any $k$ as does any three-manifold with $b_1 < 3$ by a theorem of Ozsv\'ath and Szab\'o \cite[Theorem 10.1]{OSzAnn2}.  For $\mathrm{Spin}^c$ three-manifolds $(Y,\mathfrak{s})$ with standard $HF^\infty$ there is a ``bottom-most'' correction term, denoted $d_b(Y,\mathfrak{s})$, which is defined to be the smallest grading of any non-torsion element $x \in HF^+(Y,\mathfrak{s})$ coming from an element $x' \in HF^\infty(Y,\mathfrak{s})$ which lies in the kernel of the action by $H_1(Y;\mathbb{Z})/\mathrm{Tors}$.  Notice that, in contrast to ``ordinary'' correction terms, it is not true in general that $d_b(Y,\mathfrak{s}) = -d_b(-Y,\mathfrak{s})$ (for instance, take $Y = S^2\times S^1$).  The correction terms $d_b$ give restrictions on intersection forms of negative semi-definite four-manifolds bounding a given three-manifold according to:

\begin{proposition}[Ozsv\'ath--Szab\'o \cite{OSzAbsGr}]
\label{prop:ineq} Let $Y$ be a closed oriented three-manifold (not necessarily connected) with torsion $\mathrm{Spin}^c$ structure $\mathfrak{t}$ and standard $HF^\infty$.  Then for each negative semi-definite four-manifold $W$ which bounds $Y$ so that the restriction map $H^1(W;\mathbb{Z}) \to H^1(Y;\mathbb{Z})$ is trivial, we have the inequality
\[
  c_1(\mathfrak{s})^2 + b_2^-(W) \leq 4d_b(Y,\mathfrak{t}) + 2b_1(Y)
\]for all $\mathrm{Spin}^c$ structures $\mathfrak{s}$ over $W$ whose restriction to $Y$ is $\mathfrak{t}$.
\end{proposition}

Returning to the proof of Theorem \ref{theorem:dIneq}, recall that we have a cobordism $$V_1\co S^3_1(D_-)\,\#\,M\{0,0,1\} \to S^3_1(D_+)$$ with $b_2^{\pm}(V_1) = 0$.  Therefore $c_1(\mathfrak{s})^2 = 0$ for all $\mathrm{Spin}^c$ structures $\mathfrak{s}$ on $V_1$.  Note also that $H^1(V_1;\mathbb{Z}) = 0$: consider the Mayer--Vietoris sequence applied to $W_0 = V_0\cup V_1$:
\[
  H_1(V_0\cap V_1;\mathbb{Z}) \cong \mathbb{Z}^2 \to H_1(V_0;\mathbb{Z})\oplus H_1(V_1;\mathbb{Z}) \cong \mathbb{Z}^2\oplus H_1(V_1;\mathbb{Z}) \to H_1(W_0;\mathbb{Z}) = 0
\] (recall that $V_0\cap V_1 = S^3_1(D_-)\,\#\,M\{0,0,1\}$ and that $V_0 \simeq S^3_1(D_-)\,\vee\, T^2$; $\simeq$ denoting homotopy equivalence).  Applying Proposition \ref{prop:ineq} to $V_1$ we see that:
\begin{equation}
  0 \leq 4d(S^3_1(D_+)) + 4d_b(-\left( S^3_1(D_-)\,\#\,M\{0,0,1\} \right) ) + 2\cdot2
\end{equation}

\begin{claim}\label{claim:essential} $d_b(-\left( S^3_1(D_-)\,\#\,M\{0,0,1\} \right) ) = -d(S^3_1(D_-)) + 1$.\end{claim}

Notice that this would imply Theorem \ref{theorem:dIneq}.  To prove the claim, recall that in \cite{OSzAbsGr}, Ozsv\'ath and Szab\'o calculated that $\widehat{HF}(M\{0,0,1\}) = \mathbb{Z}_{(0)}^2\oplus \mathbb{Z}_{(-1)}^2$, supported completely in the unique torsion $\mathrm{Spin}^c$ structure.  This implies that $\widehat{HF}(-M\{0,0,1\}) = \mathbb{Z}_{(0)}^2\oplus \mathbb{Z}_{(1)}^2$, and by the long exact sequence
\begin{displaymath}
\xymatrix{	\cdots \ar[r] & \widehat{HF}_i(Y,\mathfrak{s}) \ar[r] & HF_i^+(Y,\mathfrak{s}) \ar[r]^-U & HF_{i-2}^+(Y,\mathfrak{s}) \ar[r] & \cdots }
\end{displaymath}
it follows that $HF^+(-M\{0,0,1\}) = (\mathcal{T}_{(0)}^+)^2 \oplus (\mathcal{T}_{(1)}^+)^2$ where $\mathcal{T}_{(k)}$ denotes the graded $\mathbb{Z}[U]$-module $\mathbb{Z}[U,U^{-1}]/U\cdot\mathbb{Z}[U]$ graded so that multiplication by $U$ is degree $-2$ and $1\in \mathcal{T}_{(k)}$ lies in grading $k$.  Writing $HF^+(S_+^3(D_-)) = \mathcal{T}^+_{\left( d(S^3_1(D_-)) \right)} \oplus Q$ for some torsion $\mathbb{Z}[U]$ module $Q$ ($Q$ is called the \emph{reduced} Floer homology of $Y$, and is also written $HF^+_{\mathrm{red}}(Y)$), we get that $HF^+(-S^3_1(D_-)) = \mathcal{T}^+_{\left( -d(S^3_1(D_-)) \right)} \oplus Q'$, by item \ref{item:reverse} of Corollary \ref{corollary:addendum}.  By the long exact sequence
\begin{displaymath}
\xymatrix{	\cdots \ar[r] & HF_i^-(Y,\mathfrak{s}) \ar[r]& HF_i^\infty(Y,\mathfrak{s}) \ar[r] & HF^+_i(Y,\mathfrak{s}) \ar[r] & HF^-_{i-1}(Y,\mathfrak{s}) \ar[r] & \cdots }
\end{displaymath}
we get that $HF^-(-S^3_1(D_-)) = \mathcal{T}^-_{\left( -d(S^3_1(D_-))-2 \right)} \oplus Q''$ for some torsion $\mathbb{Z}[U]$-module $Q''$ where here $\mathcal{T}^-_{(k)}$ denotes the graded $\mathbb{Z}[U]$--module $U\cdot\mathbb{Z}[U]$ graded so that multiplication by $U$ is degree $-2$ and $U $ lies in grading $k$.  Using the formula
\[
	HF^-(Y_1\,\#\,Y_2,\mathfrak{s}_1\,\#\,\mathfrak{s}_2) \cong H_* \left( CF^-(Y_1,\mathfrak{s}_1) \otimes_{\mathbb{Z}[U]} CF^-(Y_2,\mathfrak{s}_2)\right)
\] from Ozsv\'ath--Szab\'o \cite{OSzAnn2}, if we use Floer homology with field coefficients $\mathbb{F}$, we have:
\begin{align*}
	HF^-(Y_1\,\#\,Y_2,\mathfrak{s}_1\,\#\,\mathfrak{s}_2) \cong HF^-(Y_1,\mathfrak{s}_1) &\otimes_{\mathbb{F}[U]} HF^-(Y_2,\mathfrak{s}_2)\\ 
&\bigoplus \mathrm{Tor}^{\mathbb{F}[U]}\left(HF^-(Y_1,\mathfrak{s}_1),HF^-(Y_2,\mathfrak{s}_2)\right)\
\end{align*} (since $\mathbb{F}[U]$ is a principal ideal domain).  It follows that 
\[
  HF^-\left( (-S^3_1(D_-))\,\#\,(-M\{0,0,1\}) \right) \cong \left( \mathcal{T}^-_{\left( -d(S^3_1(D_-))-2\right)}\right) ^2 \oplus \left( \mathcal{T}^-_{\left( -d(S^3_1(D_-))-1\right)}\right)^2 \oplus Q'''
\] for some torsion $\mathbb{F}[U]$ module $Q'''$.  Therefore
\[
	HF^+\left( (-S^3_1(D_-)) \,\#\, (-M\{0,0,1\}) \right) \cong \left( \mathcal{T}^+_{\left( -d(S^3_1(D_-))\right)}\right) ^2 \oplus \left( \mathcal{T}^+_{\left( -d(S^3_1(D_-))+1\right)}\right)^2 \oplus Q''''
\] and we have shown that 
\[
	d_b\left( -(S^3_1(D_-) \,\#\, M\{0,0,1\}) \right) = d_b\left( (-S^3_1(D_-)) \,\#\, (-M\{0,0,1\}) \right) = -d(S^3_1(D_-)) + 1
\] proving Claim \ref{claim:essential} and hence Theorem \ref{theorem:dIneq} \qedhere

\end{proof}

\section{Genus bounds}
\label{section:genus}
\begin{proof}[Proof of Theorem \ref{theorem:genus}]{\bf Step 1:} $d(S^3_1(K)) \leq 0$:

Let $g = g_4(K)$, the smooth four-ball genus of $K$, ie the minimum genus of any smooth surface smoothly embedded in the four-ball with boundary $K$.  Now attach a $-1$--framed two-handle to the four-ball along the mirror of $K$, denoted $mK$.  Now delete a small ball from the four-ball.  This gives a negative definite cobordism $S^3 \to S^3_{-1}(mK)$ whose second homology is generated by a surface of genus $g$ and square $-1$.  By item \ref{item:negDefCob} of Corollary \ref{corollary:addendum} and the fact that $d(S^3) = 0$ we get that $d(S^3_{-1}(mK)) \geq 0$.  Since $d(S^3_1(K)) = - d( -S^3_{-1}(mK))$, we are done.

{\bf Step 2:} $-d(S^3_1(K)) \leq 2g$:
 
Similar to the previous paragraph, by removing a small ball from the four-ball and then attaching a $+1$--framed two-handle to the boundary three-sphere along $K$, we obtain a cobordism
\[
	W\co S^3 \to S^3_1(K)
\] which contains a genus $g$ surface of square $+1$, $\Sigma_g$.  Let $Y_g(\pm1)$ denote an euler number $\pm1$ circle bundle over a surface of genus $g$.  In the notation of the previous section, we have $Y_1(\pm1) = M\{0,0,\pm1\}$.  $Y_g(+1)$ is of course homeomorphic to the boundary of a regular neighborhood of $\Sigma_g \subset W$.  Similar to previous discussions, by taking an internal connected sum we get a pair of cobordisms
\[
	W_1\co S^3 \to Y_g(+1)
\]and
\[
	W_2\co Y_g(+1) \to S^3_1(K)
\]with $b_2^{\pm}(W_2) = b_1(W_2) = 0$.  Notice that the oriented manifold $-Y_g(+1) = Y_g(-1)$ has standard $HF^\infty$ by Ozsv\'ath--Szab\'o \cite[Propositions 9.3 and 9.4]{OSzAbsGr} since we may connect it via a negative-definite cobordism to $\#^{2g}S^2\times S^1$, which has standard $HF^\infty$.  Applying Proposition \ref{prop:ineq}, we see that
\begin{equation}
\label{eq:first-step}
	0 \leq 4d\left(S^3_1(K)\right) + 4 d_b\left( -Y_g(+1) \right) + 2\cdot 2g.
\end{equation} Theorem \ref{theorem:genus} would follow if we could show that $d_b\left( -Y_g(+1) \right) = g$.  Indeed, we show this in Lemma \ref{lemma:computation}.  The calculation of $d_b\left( -Y_g(+1) \right)$ follows quickly from the machinery of \cite{OSint}, which we recall in Section \ref{subsection:review}.\qedhere
\end{proof}

\subsection{Review of the integer surgery formula}
\label{subsection:review}

In this section we review the essential details needed to state Ozsv\'ath and Szab\'o's ``integer surgery formula,'' referring the reader to \cite{OSint} for more details.  Suppose $(Y,\mathfrak{t})$ is a $\mathrm{Spin}^c$ three-manifold with $\mathfrak{t}$ torsion and suppose that $K\subset Y$ is a null-homologous knot.  Fixing a Seifert surface $F\subset Y$ for $K$, we can assign to $K$ its knot Floer homology $C := CFK^\infty(Y,K,F,\mathfrak{t})$, a $\mathbb{Z}\oplus\mathbb{Z}$--bifiltered chain complex well-defined up to filtered chain homotopy type as described in Ozsv\'ath--Szab\'o \cite{OSzKnot}.  This is an abelian group generated by tuples $[{\bf x},i,j]$ for integers $i,j$ and intersection points ${\bf x}$ coming from a particular Heegaard diagram for $K$ (see Ozsv\'ath--Szab\'o \cite{OSzKnot} for a proper discussion).  This group comes with an absolute $\mathbb{Q}$--grading as well as an action by $\Lambda^*H_1(Y;\mathbb{Z})/\mathrm{Tors} \otimes_\mathbb{Z} \mathbb{Z}[U,U^{-1}]$.  There is an identification of $\mathrm{Spin}^c$ structures over $Y_n(K)$ which are $\mathrm{Spin}^c$--cobordant to $\mathfrak{t}$ over a certain cobordism $W_n(K)$ with $\mathbb{Z}/n\mathbb{Z}$.  For $i\in \mathbb{Z}/n\mathbb{Z}$, let $CF^+(Y_n(K),i,\mathfrak{t})$ denote the corresponding summand of $CF^+(Y_n(K))$.

Let $A^+_{s,\mathfrak{t}} = C\{i\geq0\text{ or }j\geq s\}$ and $B_{s,\mathfrak{t}}^+ = C\{i\geq0\}$, the latter being identified with $CF^+(Y,\mathfrak{t})$.  There are maps
\[
	v^+_{s,\mathfrak{t}}\co A^+_{s,\mathfrak{t}} \to B^+_{s,\mathfrak{t}}
\]and 
\[
	h^+_{n,s,\mathfrak{t}}\co A^+_{s,\mathfrak{t}} \to B^+_{s+n,\mathfrak{t}}
\] defined as follows: $v_{s,\mathfrak{t}}$ is just the projection $C\{i\geq0\text{ or }j\geq s\} \to C\{i\geq0\}$ while $h^+_{s,\mathfrak{t}}$ is the projection $C\{i\geq0\text{ or }j\geq s\} \to C\{j\geq s\}$ followed by an identification $C\{j\geq s\} \cong C\{j\geq0\}$ (induced by multiplication by $U^s$) followed by a ``natural'' homotopy equivalence $h\co C\{j\geq0\} \to C\{i\geq0\}$.  The map $h$ is obtained by the handleslide invariance of Heegaard Floer homology and is natural in the sense that the induced map on homology is independent (up to a sign) of a chosen sequence of handleslides.  Set
\[
	\mathbb{A}^+_{i,\mathfrak{t}} = \bigoplus_{\{s\in\mathbb{Z}\,|\,s\equiv i\mod n\}} A^+_{s,\mathfrak{t}}
\]and:
\[
	\mathbb{B}^+_{i,\mathfrak{t}} = \bigoplus_{\{s\in\mathbb{Z}\,|\,s\equiv i\mod n\}} B^+_{s,\mathfrak{t}}
\]  Define 
\[
	\mathcal{D}^+_{i,\mathfrak{t},n}\co \mathbb{A}^+_{i,\mathfrak{t}} \to \mathbb{B}^+_{i,\mathfrak{t}}
\] by 
\[
	 \mathcal{D}^+_{i,\mathfrak{t},n}\{a_s\}_{s\in\mathbb{Z}} = \{ v^+_{s,\mathfrak{t}}(a_s) + h^+_{s-n,\mathfrak{t}}(a_{s-n})\}_{s\in\mathbb{Z}}.
\]  Assign gradings to $\mathbb{A}^+_{i,\mathfrak{t}}$, $\mathbb{B}^+_{i,\mathfrak{t}}$ as follows.  Under the identification $B^+_{s,\mathfrak{t}} \cong CF^+(Y,\mathfrak{t})$, we map homogeneous elements of degree $d$ in $CF^+(Y,\mathfrak{t})$ to homogeneous elements of $B^+_{s,\mathfrak{t}}$ of degree
\begin{equation}
\label{equation:gr1} 
d+2\sigma l + n\ell(\ell-1)-1, \;\mathrm{where}\;\; 0 \leq \sigma < n \;\;\mathrm{and}\;\; s = \sigma + n\ell  \;\;\mathrm{if}\;\;  n>0,
\end{equation} or
\begin{equation}
\label{equation:gr2} d - 2\sigma\ell +n\ell(\ell-1) \;\;\mathrm{where}\;\; 0 \leq \sigma < -n \;\;\mathrm{and}\;\; s = -(\sigma - n\ell) \;\;\mathrm{if}\;\; n<0.
\end{equation}
It is then possible to assign gradings to the $A^+_{s,\mathfrak{t}}$ which are consistent with their natural relative $\mathbb{Z}$--gradings in such a way that the maps $v^+_{s,\mathfrak{t}}$ and $h^+_{s,\mathfrak{t}}$ are homogeneous of degree $-1$.  With this all in place, we may now state the ``integer surgery formula'' of \cite{OSzAbsGr}:
\begin{theorem}[Ozsv\'ath--Szab\'o \cite{OSzAbsGr}]
\label{theorem:intSurgery}
Fix a $\mathrm{Spin}^c$ structure $\mathfrak{t}$ over $Y$ whose first Chern class is torsion, $K\subset Y$ a null-homologous knot, and $n$ a non-zero integer.  For each $i \in \mathbb{Z}/n\mathbb{Z}$, the mapping cone $\mathbb{X}^+_{i,\mathfrak{t}}(n)$ of \[\mathcal{D}^+_{n,i,\mathfrak{t}}\co \mathbb{A}^+_{i,\mathfrak{t}} \to \mathbb{B}^+_{i,\mathfrak{t}}\] is isomorphic, as a relatively graded  $\Lambda^*H_1(Y;\mathbb{Z})/\mathrm{Tors} \otimes_{\mathbb{Z}} \mathbb{Z}[U]$--module, to $CF^+(Y_n(K),i,\mathfrak{t})$.  In fact, this isomorphism $\mathbb{X}^+_{i,\mathfrak{t}(n)} \to CF^+(Y_n(K),i,\mathfrak{t})$ is homogeneous of degree $d(n,i)$ where 
\[
	d(n,i) = -\max_{\{s\in\mathbb{Z}\,|\,s\equiv i\mod n\}} \frac{1}{4}\left( 1 - \left( \frac{n+2s}{n} \right)^2 \right)
\] for $n>0$ and $d(n,i) = -d(-n,i)$ for $n<0$.
\end{theorem}
Recall that the mapping cone of the map $\mathcal{D}^+_{n,i,\mathfrak{t}}\co \mathbb{A}^+_{i,\mathfrak{t}} \to \mathbb{B}^+_{i,\mathfrak{t}}$ has underlying group $\mathbb{X}^+_{i,\mathfrak{t}}(n) = \mathbb{A}^+_{i,\mathfrak{t}} \oplus \mathbb{B}^+_{i,\mathfrak{t}}$ and differential
\[
	\partial_{\mathbb{X}^+_{i,\mathfrak{t}}(n)} = \begin{pmatrix} \partial_{\mathbb{A}^+_{i,\mathfrak{t}}} & 0 \\
								\mathcal{D}^+_{n,i,\mathfrak{t}} & \partial_{\mathbb{B}^+_{i,\mathfrak{t}}} \end{pmatrix}.
\]  Also, when $n=\pm1$, $\mathrm{Spin}^c(Y_n(K)) \cong \mathrm{Spin}^c(Y)$ and there is no additional choice of $i\in\mathbb{Z}/n\mathbb{Z}$.  When this is satisfied, we write simply $\mathbb{X}^+_{\mathfrak{t}}(n)$ instead of $\mathbb{X}^+_{i,\mathfrak{t}}(n)$.  

\subsection{A useful computation}
\label{section:useful_computation}

\begin{lemma}
\label{lemma:computation} 
For $Y_g(\pm1)$ as before, we have $$d_b(Y_g(\pm1);\mathbb{Z}_2) = \mp g$$ where $d(Y,\mathfrak{s};\mathbb{Z}_2)$ denotes the $d$--invariant of $(Y,\mathfrak{s})$ as computed from Floer homology with coefficients in $\mathbb{Z}_2$.     
\end{lemma}
\begin{proof}
The oriented manifold $-Y_g(+1) = Y_g(-1)$ may be obtained as $-1$--surgery on the knot $K$ (the ``Borromean knot'') shown in Figure \ref{figure:borromean}.

\begin{figure}
\labellist
\small\hair 2pt
\pinlabel $0$ at 50 90
\pinlabel $0$ at 152 90
\pinlabel $0$ at 76 70
\pinlabel $g$ at 101 19
\pinlabel $0$ at 183 64
\pinlabel $K$ at 101 -6
\endlabellist
\begin{center}
\includegraphics{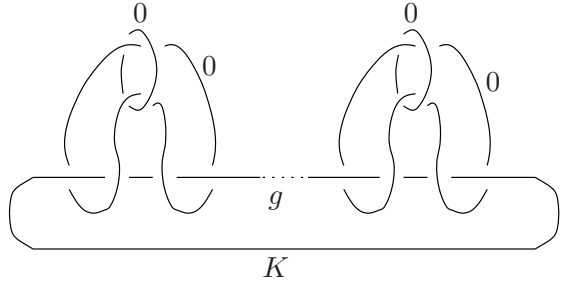}
\caption{\label{figure:borromean} The Borromean knot $K$.    }
\end{center}
\end{figure}

We start with the calculation for $Y_g(-1)$.  Since the second homology of $Y_g(-1)$ is generated by embedded tori, the adjunction inequality (Ozsv\'ath--Szab\'o \cite[Theorem 7.1]{OSzAnn2}) implies that $HF^+(Y_g(+1),\mathfrak{t})$ is non-zero only in the unique torsion $\mathrm{Spin}^c$ structure $\mathfrak{t}_0\in \mathrm{Spin}^c(Y_g(-1))$.  The knot Floer complex for the Borromean knot $K$ is calculated in Ozsv\'ath--Szab\'o \cite{OSzKnot} to be 
\[
	C := CFK^\infty(\#^{2g}S^2\times S^1,K) \cong \Lambda^*H^1(\Sigma_g;\mathbb{Z}) \otimes_\mathbb{Z} \mathbb{Z}[U,U^{-1}]
\] with $\mathbb{Z}\oplus\mathbb{Z}$--bifiltration given by:
\begin{equation}
\label{equation:filtration}
	C\{i,j\} = \Lambda^{g-i+j}H^1(\Sigma_g;\mathbb{Z})\otimes_\mathbb{Z}U^{-i}
\end{equation}  Furthermore, the group $C\{i,j\}$ is supported in grading $i+j$ and all differentials vanish (including all ``higher'' differentials coming from the spectral sequence $HFK^\infty \Rightarrow HF^\infty$).  

Under the above identification, and the identification $H_1(\Sigma_g;\mathbb{Z}) \cong H_1(\#^{2g}S^2\times S^1;\mathbb{Z})$, the action of $\gamma\in H_1(\Sigma_g;\mathbb{Z})$ on $C \cong \Lambda^*H^1(\Sigma_g;\mathbb{Z}) \otimes_\mathbb{Z} \mathbb{Z}[U,U^{-1}]$ is given explicitly by
\begin{equation}
\label{equation:action}
	\gamma\cdot(\omega\otimes U^j) = \iota_\gamma\omega\otimes U^j + \mathrm{PD}(\gamma)\wedge\omega\otimes U^{j+1}
\end{equation} where $\iota_\gamma$ denotes contraction.  Since $HF^\infty(\#^{2g}S^2\times S^1) \cong C$, this action may be viewed as a reflection of the fact that $HF^\infty(\#^{2g}S^2\times S^1)$ is standard.

The only presumably non-combinatorial ingredient in the integer surgery formula (once the complex $C = CFK^\infty(Y,K,i)$ is at hand) is the necessary explicit identification of the natural homotopy equivalence $h\co C\{j\geq0\} \to C\{i\geq0\}$.  The homotopy $h$ takes a particularly simple form for Floer homology with coefficients in $\mathbb{Z}_2$, with which we work for the remainder of this section.  The description of $h$ is as follows:
\begin{proposition}[Ozsv\'ath--Szab\'o \cite{OSzAbsGr}]
\label{prop:htpy}
For the Borromean knot $K$, the natural homotopy equivalence $h\co C\{j\geq0\} \to C\{i\geq0\}$ sends $C\{j,i\}$ to $C\{i,j\}$.  
\end{proposition} Interestingly, the above proposition does not hold for Floer homology with coefficients in $\mathbb{Z}$ (see Jabuka--Mark \cite{JM2} for a description).   

We picture the complex $\mathbb{X}^+(-1)$ as below:
\begin{displaymath}
\xymatrix@R=4pc@C=4pc{	
  A^+_{-1} \ar[d]_{\cdots\;\;\;v^+_{-1}}^{} & A^+_0 \ar[dl]_{h^+_0} \ar[d]_-{v^+_0} & A^+_1 \ar[dl]_{h^+_1} \ar[d]_-{v^+_1}^{\;\;\; \cdots} \\
  B^+_{-1}				& B^+_0					& B^+_1
}\end{displaymath} For simplicity of discussion, we currently restrict to the case of $g=1$.  In this case, a piece of $\mathbb{X}^+(-1)$ looks like Figure \ref{figure:wedges}.

\begin{figure}[ht!]
\labellist
\small\hair 2pt


\pinlabel $B^+_{-1}$ at 10 70
\pinlabel $\Lambda^2$ at 36 52
\pinlabel $\Lambda^1$ at 36 37
\pinlabel $\Lambda^0$ at 36 22
\pinlabel $\Lambda^2$ at 51 67
\pinlabel $\Lambda^1$ at 51 52
\pinlabel $\Lambda^0$ at 51 37
\pinlabel $\iddots$ at 66 67

\pinlabel $B^+_0$ at 100 70
\pinlabel $\Lambda^2$ at 126 52
\pinlabel $\Lambda^1$ at 126 37
\pinlabel $\Lambda^0$ at 126 22
\pinlabel $\Lambda^2$ at 141 67
\pinlabel $\Lambda^1$ at 141 52
\pinlabel $\Lambda^0$ at 141 37
\pinlabel $\iddots$ at 156 67

\pinlabel $B^+_1$ at 190 70
\pinlabel $\Lambda^2$ at 216 52
\pinlabel $\Lambda^1$ at 216 37
\pinlabel $\Lambda^0$ at 216 22
\pinlabel $\Lambda^2$ at 231 67
\pinlabel $\Lambda^1$ at 231 52
\pinlabel $\Lambda^0$ at 231 37
\pinlabel $\iddots$ at 246 67


\pinlabel $A^+_{-1}$ at 10 178
\pinlabel $\Lambda^2$ at 6 130
\pinlabel $\Lambda^1$ at 21 130
\pinlabel $\Lambda^2$ at 21 145
\pinlabel $\Lambda^2$ at 36 160
\pinlabel $\Lambda^1$ at 36 145
\pinlabel $\Lambda^0$ at 36 130
\pinlabel $\Lambda^2$ at 51 175
\pinlabel $\Lambda^1$ at 51 160
\pinlabel $\Lambda^0$ at 51 145
\pinlabel $\iddots$ at 66 175

\pinlabel $A^+_0$ at 100 178
\pinlabel $\Lambda^2$ at 111 145
\pinlabel $\Lambda^2$ at 126 160
\pinlabel $\Lambda^1$ at 126 145
\pinlabel $\Lambda^0$ at 126 130
\pinlabel $\Lambda^2$ at 141 175
\pinlabel $\Lambda^1$ at 141 160
\pinlabel $\Lambda^0$ at 141 145
\pinlabel $\iddots$ at 156 175

\pinlabel $A^+_1$ at 190 178
\pinlabel $\Lambda^2$ at 216 160
\pinlabel $\Lambda^1$ at 216 145
\pinlabel $\Lambda^0$ at 216 130
\pinlabel $\Lambda^2$ at 231 175
\pinlabel $\Lambda^1$ at 231 160
\pinlabel $\Lambda^0$ at 231 145
\pinlabel $\iddots$ at 246 175

\endlabellist
\begin{center}
\includegraphics{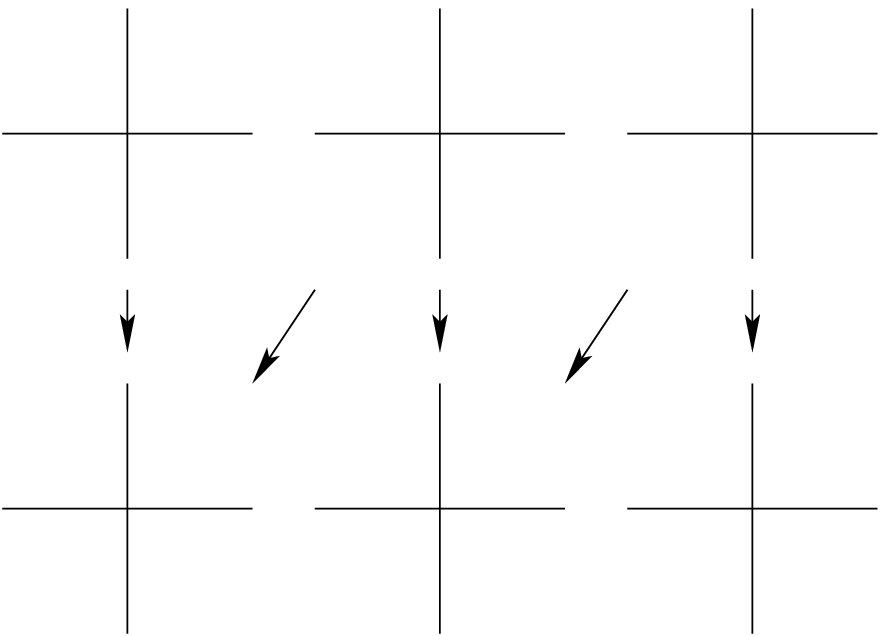}
\caption{ \label{figure:wedges}A portion of the complex $\mathbb{X}(-1)$.  We suppress the $U$'s from the notation, since they can be determined from the position in the plane, according to Equation \ref{equation:filtration}. }
\end{center}
\end{figure}

We claim that the correction terms of $Y_1(-1)$ can be read off from Figure \ref{figure:wedges}.  Indeed, writing $H_1(\Sigma_1;\mathbb{Z}) = \mathbb{Z}\langle a,b\rangle$ and $\Lambda^*H_1(\Sigma_1;\mathbb{Z}) = \mathbb{Z}\langle1,a,b,a\wedge b\rangle$, consider the element $1\otimes U^{-1} \in B^+_0\{1,0\}$.  Since 
$$
	\partial_{\mathbb{X}^+(-1)} = \begin{pmatrix} \partial_{\mathbb{A}^+} & 0\\
						\mathcal{D}^+ & \partial_{\mathbb{B}^+}\end{pmatrix} = \begin{pmatrix} 0 & 0\\
												\mathcal{D}^+ & 0\end{pmatrix},
$$ it follows that $1\otimes U^{-1} \in B^+_0$ is a cycle.  We claim that $1\otimes U^{-1}\in B^+_0$ is also not a boundary.  Indeed, suppose that $\partial_{\mathbb{X}^+(-1)}(x) = 1\otimes U^{-1} \in B^+_0$ for some $x\in\mathbb{X}^+(-1)$.  Then $x$ would necessarily have either a non-zero component in $A^+_0\{1,0\}$ or a non-zero component in $A^+_1\{1,2\}$.  In either case, a simple diagram chase shows that $x$ cannot be extended to a cycle in $\mathbb{X}^+(-1)$ (ride the zig-zag and notice that $x$ would have infinitely many non-zero components in $\mathbb{A}^+$).

It is, however, the case that $U\cdot(1\otimes U^{-1}) = 1\otimes U^0 \in B^+_0\{0,-1\}$ is a boundary: the element $1\otimes U^0 \in A^+_0\{0,-1\}$ maps to it under $\partial_{\mathbb{X}^+(-1)}$.

Using similar reasoning, one can show that the elements 
\[
	1\otimes U^{-1}\in B^+_0\{1,0\},\; a\otimes U^0\in B^+_0\{0,0\},\; b\otimes U^0\in B^+_0\{0,0\},\; \;\mathrm{and}\; (a\wedge b)\otimes U^0 \in B^+_0\{0,1\}
\] all represent generators for the four ``towers'' of $HF^+(Y_1(-1))$.  According to the grading formula, Equation \ref{equation:gr2}, these elements have grading 1,0,0, and 1, respectfully.  It follows that:
\[
	HF^+(Y_{1}(-1)) = \left( \mathcal{T}_{(0)}^+ \right)^2 \oplus \left( \mathcal{T}^+_{(+1)} \right)^2 \oplus HF^+_{\mathrm{red}}(Y_{1}(-1))
\] where here we are using a slight abuse of notation: $\mathcal{T}^+$ now denotes $\mathcal{T}^+\otimes\mathbb{Z}_2$ (for the previous definition of $\mathcal{T}^+$) and Floer homology with $\mathbb{Z}_2$--coefficients is understood.   Further, according to the action, Equation \ref{equation:action}, it follows that $d_b(Y_1(-1);\mathbb{Z}_2) = 1$ (alternatively this follows since $Y_1(-1)$ has standard $HF^\infty$ via Equation \ref{equation:std}).  Although we already knew how to compute this, the advantage of this calculation is that the reasoning generalizes to arbitrary $g$.  Indeed, for general $g$ one can check that the intersection
\[
	B^+_0\{\left(i\geq0\;\mathrm{and}\;j=0\right) \;\mathrm{or}\; \left(i=0 \;\mathrm{and}\; j\geq0\right)\}
\] gives representatives for generators of the $2^{2g}$ ``towers'' of $HF^+(Y_g(-1))$.  Using the grading formula, Equation \ref{equation:gr2}, it follows that:
\begin{align*}
	HF^+(Y_{g}(-1))  \cong \left( \mathcal{T}_{(0)}^+ \right)^{ \bigl( \begin{smallmatrix}2g\\g\end{smallmatrix} \bigr) } &\oplus \left( \mathcal{T}_{(1)}^+ \right)^{ 2\bigl( \begin{smallmatrix}2g\\g-1\end{smallmatrix}\bigr) } \oplus \left( \mathcal{T}_{(2)}^+ \right)^{ 2\bigl(\begin{smallmatrix}2g\\g-2\end{smallmatrix}\bigr) } \oplus \cdots\\ \cdots&\oplus \left( \mathcal{T}_{(g-1)}^+ \right)^{ 2\bigl(\begin{smallmatrix}2g\\1\end{smallmatrix}\bigr) } \oplus \left( \mathcal{T}_{(g)}^+ \right)^{ 2 } \oplus HF^+_{\mathrm{red}}(Y_g(-1))
\end{align*}
By the action formula, Equation \ref{equation:action}, (or the fact that $HF^\infty(Y_g(-1))$ is standard) it follows that \[ d_b(Y_g(-1);\mathbb{Z}_2) = g.\]  Calculating $d_b(Y_g(-1);\mathbb{Z}_2)$ is similar.  In that case, the generators of the intersection
\[
	A^+_0\{\left(i\leq0\;\mathrm{and}\;j=0\right) \;\mathrm{or}\; \left(i=0 \;\mathrm{and}\; j\leq0\right)\}
\] may be extended (in one step) to representative cycles for the homology $HF^+(Y_g(+1))$.  Using Equation \ref{equation:gr1} one calculates that:
\begin{align*}
	HF^+(Y_{g}(+1))  \cong \left( \mathcal{T}_{(0)}^+ \right)^{ \bigl( \begin{smallmatrix}2g\\g\end{smallmatrix} \bigr) } &\oplus \left( \mathcal{T}_{(-1)}^+ \right)^{ 2\bigl( \begin{smallmatrix}2g\\g-1\end{smallmatrix}\bigr) } \oplus \left( \mathcal{T}_{(-2)}^+ \right)^{ 2\bigl(\begin{smallmatrix}2g\\g-2\end{smallmatrix}\bigr) } \oplus \cdots\\ \cdots&\oplus \left( \mathcal{T}_{(-g+1)}^+ \right)^{ 2\bigl(\begin{smallmatrix}2g\\1\end{smallmatrix}\bigr) } \oplus \left( \mathcal{T}_{(-g)}^+ \right)^{ 2 } \oplus HF^+_{\mathrm{red}}(Y_g(+1))
\end{align*}
The action formula, Equation \ref{equation:action}, gives 
\begin{equation}
\label{eq:dYg}
 d_b(Y_g(+1);\mathbb{Z}_2) = -g
\end{equation}
proving Lemma \ref{lemma:computation}\qedhere
\end{proof}


An alternative approach to the calculation of the correction terms of $\pm1$--surgery on the Borromean knot of genus $g$ is the integer surgery exact sequence, together with Jabuka-Mark's calculation of the Floer homology of $S^1\times\Sigma_g$ \cite{JM3}.  It is interesting to note that $HF^+_\mathrm{red}(Y_g(-1))$ has been calculated in Ozsv\'ath--Szab\'o \cite{OSzKnot}.

In fact, the methods in this section give the following:
\begin{theorem}
\label{theorem:b2plus}
Let $X$ be a smooth, simply-connected, compact, oriented four-manifold with a homology sphere $Y$ as boundary with $b_2^+(X)=1$ and $b_2^-(X)=0$.  Let $\Sigma\subset X$ be a closed surface of genus $g$ and self-intersection $\Sigma\cdot\Sigma=1$.  Then $0\leq d(Y;\mathbb{Z}/2)+2g.$
\end{theorem}

\section{Computations}
\label{section:computations}
In this section we discuss an algorithm to compute the invariants $d(S^3_{\pm1}(K))$ assuming we know the filtered chain homotopy type of the knot complex $CFK^\infty(K)$ of Ozsv\'ath--Szab\'o \cite{OSzKnot}.  We also discuss a computer implementation of this algorithm.  The algorithm we use is based on the theory of Ozsv\'ath--Szab\'o \cite{OSzKnot,OSzAbsGr}, and Rasmussen \cite{Ra} and has three steps:
\begin{enumerate}
\item \label{item:MAIN} Use the theory of Ozsv\'ath--Szab\'o \cite{OSzKnot} or Rasmussen \cite{Ra} to compute the graded $\mathbb{F}[U]$-module \[HF^+(S^3_{-p}(K),\mathfrak{s}_0)\] where $-p$ is ``sufficiently negative'' and $\mathfrak{s}_0$ is a suitable $\mathrm{Spin}^c$ structure on $S^3_{-p}(K)$ (and similarly for $S^3_{p}(K)$ for $p$ ``sufficiently positive'').
\item \label{item:EASY1} Use exact sequences to compute the correction terms of $S^3_0(K)$.
\item \label{item:EASY2} Use a simple relation between the correction terms of $S^3_0(K)$ and the correction terms of $S^3_{\pm1}(K)$.
\end{enumerate}

We describe the steps in reverse order.  Recall from \cite{OSzAbsGr} that for a closed oriented three-manifold $Y_0$ with $H_1(Y_0;\mathbb{Z})\cong\mathbb{Z}$, there are two correction terms $d_{\pm1/2}(Y_0)$, where $d_{\pm1/2}(Y_0)$ is the minimal grading of any non-torsion element in the image of $HF^\infty(Y_0,\mathfrak{t}_0)$ in $HF^{\pm}(Y_0,\mathfrak{t}_0)$ with grading $\pm1/2$ modulo 2.  Then step \ref{item:EASY2} follows from Ozsv\'ath--Szab\'o \cite[Proposition 4.12 ]{OSzAbsGr} which states that 
\begin{equation}
\label{equation:d1/2}
	d_{1/2}(S^3_0(K)) - 1/2 = d(S^3_1(K))
\end{equation} and 
\begin{equation}
\label{equation:d-1/2}
	d_{-1/2}(S^3_0(K)) + 1/2 = d(S^3_{-1}(K)).
\end{equation}This proposition is an easy consequence of the fact that $HF^\infty(S^3_0(K))$ is standard and the exact sequence, Equation \ref{equation:exact}.

For step \ref{item:EASY1}, recall the integral surgeries long exact sequence (see Ozsv\'ath--Szab\'o \cite[Theorem 9.19]{OSzAnn2}; \cite{OSzAbsGr} for the graded version).  Let $K\subset Y$ be a knot in an integral homology three-sphere and $p$ a positive integer.  Then we get a map 
\[
	Q\co \mathrm{Spin}^c(Y_0) \to \mathrm{Spin}^c(Y_p)
\] and a long exact sequence of the form
\begin{equation}
\label{equation:pexact}
\xymatrix{
	\cdots \ar[r]^-{F_1} &HF^+(Y_0,[\mathfrak{t}]) \ar[r]^{F_2} &HF^+(Y_p,\mathfrak{t}) \ar[r]^{F_3} &HF^+(Y) \ar[r] &\cdots}
\end{equation} where 
\[
	HF^+(Y_0,[\mathfrak{t}]) = \bigoplus_{\mathfrak{t}'\in Q^{-1}(\mathfrak{t})} HF^+(Y_0,\mathfrak{t}').
\]  Moreover, the component of $F_1$ in the above exact sequence which takes $HF^+(Y)$ into the $\mathfrak{t}_0$-component of $HF^+(Y_0,[Q(\mathfrak{t}_0)])$ has degree $-1/2$, while the restriction of $F_2$ to the $HF^+(Y_0,\mathfrak{t}_0)$--summand of $HF^+(Y_0,[Q(\mathfrak{t}_0)])$ has degree $\left(\frac{p-3}{4}\right)$.  It follows from this exact sequence that $d_{1/2}(S^3_0(K)) + \left(\frac{p-3}{4}\right) = d(S^3_{+p}(K),\mathfrak{s}_0)$.  To get $d_{-1/2}(S^3_0(K))$, we may use a similar exact sequence for negative surgeries.

 $CFK^\infty(K)$ is finitely generated as a complex over $\mathbb{Z}[U,U^{-1}]$.  It comes with an absolute $\mathbb{Z}$--grading and a $\mathbb{Z}\oplus\mathbb{Z}$--bifiltration.  Generators are written $[{\bf x},i,j]$ for integers $i,j$ and intersection points ${\bf x}$ in a Heegaard diagram for $K$.  The $U$--action is given by $U\cdot[{\bf x},i,j] = [{\bf x},i-1,j-1]$.  We picture these complexes as graphs in the plane, as in Figure \ref{figure:t34} (1).  Dots represent generating $\mathbb{Z}$'s while arrows represent differentials.  The absolute grading is (basically) pinned down by the fact that if we consider the ``$y$--slice'' quotient complex $CFK^\infty\{i=0\}$, then its homology (which is guaranteed to be a single copy of $\mathbb{Z}$---the generator of $\widehat{HF}(S^3)\cong\mathbb{Z}_{(0)}$) is supported in grading $0$, the fact that the $U$--action drops absolute grading by $2$, and the fact that differentials drop grading by $1$\footnote{This is actually not quite true: not all vertices appearing in the complex $CFK^\infty$ are related to the generator of $\widehat{HF}(S^3)$ by a sequence of $U$--maps and differentials.  To grade these remaining vertices, we have to go back to the Heegaard diagram.  However, for the purpose of computing $d$--invariants, we do not need to look at these at all.\label{fn:traversal}}.  

In order to accomplish step \ref{item:MAIN}, we use the following theorem of Ozsv\'ath--Szab\'o, Rasmussen (see, for instance, Ozsv\'ath--Szab\'o \cite[Corollary 4.3]{OSzKnot})
\begin{proposition}[Ozsv\'ath--Szab\'o \cite{OSzKnot}, Rasmussen \cite{Ra}]
\label{prop:largeN}
 Let $K$ be a knot in the three-sphere.  Then there exists a positive integer $N$ with the property that for all $p\geq N$ we have that 
\[
	HF^+_\ell(S^3_{-p}(K),[0]) \cong H_k(CFK^\infty(K)\{i\geq0\text{ and }j\geq0\}
\] where 
\[
	\ell = k + \left( \frac{1-p}{4} \right).
\]Similarly,
\[
	HF^+_\ell(S^3_{p}(K),[0]) \cong H_k(CFK^\infty(K)\{i\geq0\text{ or }j\geq0\})
\] where 
\[
	\ell = k + \left( \frac{p-1}{4} \right).
\]In fact, we can take $N=2g-1$ where $g = g(K)$ is the knot genus.
\end{proposition}
We can pick off the correction terms $d(S^3_{\pm p}(K),0)$ from this theorem:  we know that \[H_*(CFK^\infty(K))\cong \mathbb{Z}[U,U^{-1}]\] so choose a class $a\in CFK^\infty(K)$ which generates this homology as a $\mathbb{Z}[U,U^{-1}]$--module.  Look at a sufficiently negative $U$--power of this generator.  This generates the ``tower'' of $H_*(CFK^\infty(K)\{i\geq0\text{ or }j\geq0\}$.  All one has to do is start taking $U$--powers of $[a] \in H_k(CFK^\infty(K)\{i\geq0\text{ or }j\geq0\}$ and see when they vanish in homology.  The grading of the last-surviving $U$--power of $[a]$ is $d(S^3_{+p}(K),[0])$ (after the grading shift).  A similar story allows one to compute $d(S^3_{-p}(K),0)$.

Write $\widetilde{d}(S^3_{p}(K),[0]) = d(S^3_p(K),[0]) - \frac{p-1}{4}$, the ``unshifted'' correction term of the group $H_k(CFK^\infty(K)\{i\geq0\text{ or }j\geq0\})$.  Putting together the previous discussion, by the exact sequence, Equation \ref{equation:pexact} we have:
\[
	d_{\frac{1}{2}}(S^3_0(K)) + \frac{p-3}{4} = d(S^3_p(K),[0]) = \widetilde{d}(S^3_{p}(K),0) + \frac{p-1}{4}
\] Using Equation \ref{equation:d1/2} we get:
\[
	\frac{1}{2} + d(S^3_1(K)) + \frac{p-3}{4} = \widetilde{d}(S^3_{p}(K),0) + \frac{p-1}{4}
\] ie
\begin{equation}
\label{equation:shift}
	d(S^3_1(K)) = \widetilde{d}(S^3_{p}(K),[0]).
\end{equation}

We now discuss how to teach a computer to do step \ref{item:MAIN}.  In fact, we implemented this in \CC\; in a program called \texttt{dCalc} (\texttt{beta}).  The source code is available at \[\texttt{http://www.math.columbia.edu/\textasciitilde tpeters}.\]  As previously mentioned, $CFK^\infty(K)$ is finitely generated as a complex over $\mathbb{Z}[U,U^{-1}]$.  By a symmetry property of the knot Floer homology, we may assume that the corresponding graph is symmetric about the line $i=j$.  With such a graph at hand, we choose a generating set which is
\begin{enumerate}
\item Minimal: no smaller subset of it generates $CFK^\infty(Y)$.
\item In the first quadrant, $i\geq0$ and $j\geq0$.
\item As close to the origin as possible.
\end{enumerate} (See Figure \ref{figure:t34} for an example). 
We started with a digraph data structure to represent the complexes.  Vertices were marked with $\mathbb{Z}\oplus\mathbb{Z}$--bifiltration levels and could be marked with gradings.  Vertices are also marked to keep track of bases.  The first step was to fill in the gradings.  For this we needed to compute the ``$y$--slice'' described before.  This is determined by a finite number of $U$--translates of our chosen generating set.  Once we find the generator of $\widehat{HF}(S^3)$, it is a problem in graph traversal to fill in (most of---see footnote \ref{fn:traversal}) the other gradings.  Our chosen generating set will necessarily have $1$--dimensional homology (over $\mathbb{Z}_2$) and its generator $x$ will have a grading computed from the graph traversal.  This generator maps to a generator of $HF^+(S^3) \cong \mathcal{T}^+_{(0)}$.  To compute $d(S^3_{+1}(K))$, we start by taking a finite piece of the complex $CFK^\infty(K)\{i\geq0\text{ or }j\geq0\}$ (more specifically, take our chosen generating set and start hitting it by $U$---at some point it will disappear out of the ``hook'' region $\{i\geq0\text{ or }j\geq0\}$.  Take only those images which appear in the hook).  In the graph implementation, this just involves shifting filtration levels and throwing away some vertices as they exit the hook.  Now take the generator $x$ and start pushing it down by $U$ until it dies in homology.  Its grading just before it dies is $d(S^3_{+1}(K))$, by Equation \ref{equation:shift}.  To compute $d(S^3_{-1}(K))$, one runs a similar story, but instead of using the ``hook region'' $\{i\geq0\text{ or }j\geq0\}$, one uses the first quadrant $\{i\geq0\text{ and }j\geq0\}$.

Since one knows how knot complexes behave under connected sum of knots (tensor product over $\mathbb{Z}[U,U^{-1}]$; see Ozsv\'ath--Szab\'o \cite[Theorem 7.1]{OSzKnot} for the precise formulation), we implemented this as well, allowing users to compute correction terms of surgeries on connected sums of knots.

\subsection{A few examples}
\label{section:examples}
\begin{figure}
\labellist
\small\hair 2pt
\pinlabel $\mbox{1. Find minimal generating.}$ at 65 308
\pinlabel $\mbox{  graph (shown in bold).}$ at 65 296
\pinlabel $\mbox{2. Find generator of ``$y$-slice.''}$ at 236 306
\pinlabel $\mbox{Here it is!}$ at 251 413
\pinlabel $\mbox{3. Fill in gradings.}$ at 45 145
\pinlabel $\mbox{This will do.}$ at 266 188
\pinlabel $\mbox{4. Find a generator of}$ at 216 150
\pinlabel $\mbox{   homology.}$ at 216 139
\pinlabel $\mbox{Not quite}$ at 106 30
\pinlabel $\mbox{dead yet...}$ at 122 18
\pinlabel $\mbox{5. Start pushing generator down}$ at 80 -6
\pinlabel $\mbox{   by $U$ until it dies in the}$ at 70 -18
\pinlabel $\mbox{   ``hook'' region, $i,j\geq0$.}$ at 65 -30
\pinlabel $\mbox{RIP}$ at 248 18
\pinlabel $\mbox{6. Grading before death was}$ at 236 -6
\pinlabel $-2\mbox{ so }d(S^3_1(T_{3,4}))=-2.$ at 236 -18
\pinlabel $i$ at 145 380
\pinlabel $j$ at 66 450
\pinlabel $0$ at 52 273
\pinlabel $1$ at 92 278
\pinlabel $0$ at 74 234
\pinlabel $1$ at 122 248
\pinlabel $0$ at 113 213
\endlabellist
\begin{center}
\includegraphics{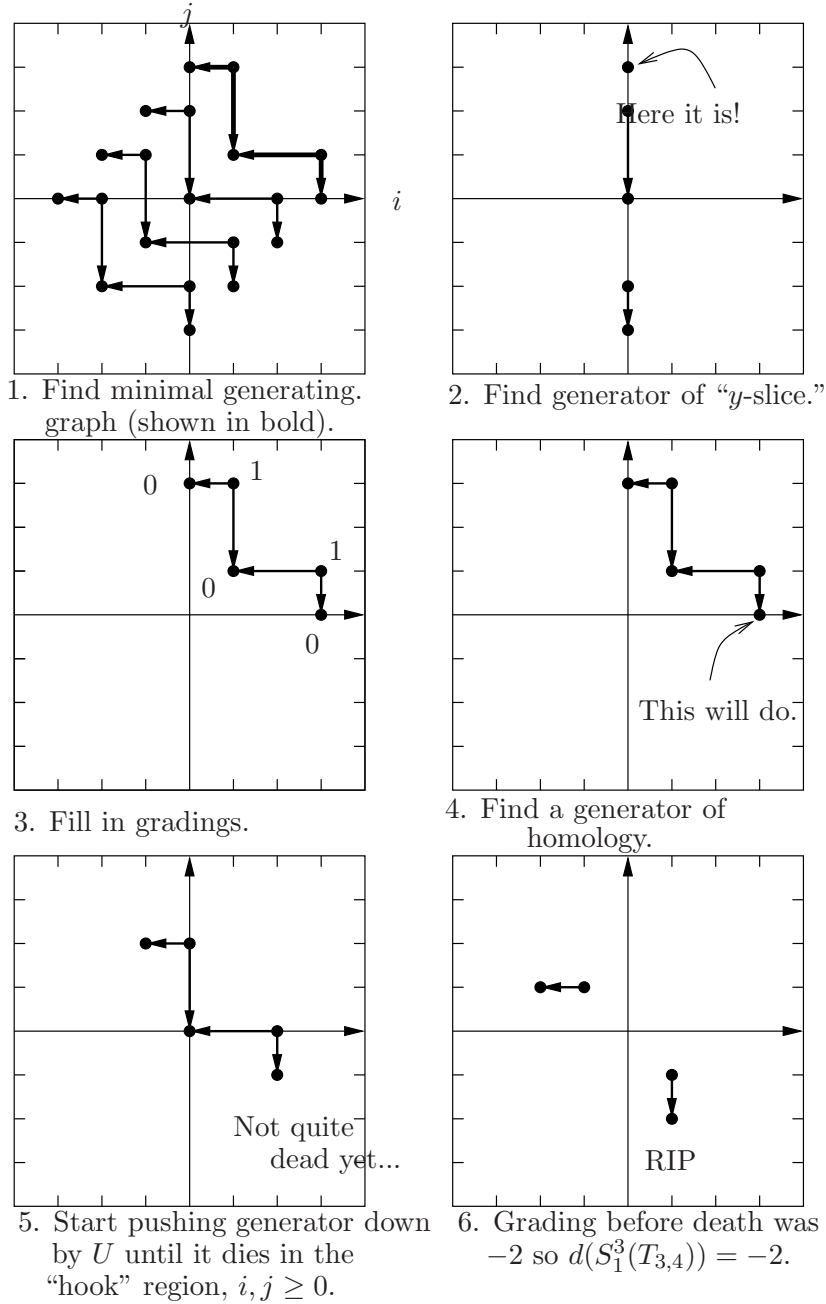}
\caption{\label{figure:t34} An example of the algorithm used, applied to the $(3,4)$--torus knot, $T_{3,4}$.    }
\end{center}
\end{figure}

In this section, Floer homology with mod-2 coefficients is understood.  Figure \ref{figure:t34} shows an example of the algorithm described in Section \ref{section:computations}.  Here, we are given the knot complex $CFK^\infty(T_{3,4})$ for the $(3,4)$--torus knot, which was computed in Ozsv\'ath--Szab\'o \cite[Section 5.1]{OSzAbsGr} or by \cite[Theorem 1.2]{OSLens}.

\begin{figure}
\labellist
\small\hair 2pt
\pinlabel $\mathbb{Z}_2$ at 30 39
\pinlabel $\mathbb{Z}_2$ at 8 19
\pinlabel $\mathbb{Z}_2$ at 48 55
\pinlabel Maslov at 25 -3 
\pinlabel Alexander at -35 55
\endlabellist
\begin{center}
\includegraphics{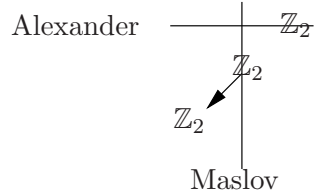}
\caption{\label{figure:tref}The $E^1$ page of the spectral sequence $\widehat{HFK}(RHT)\Rightarrow\widehat{HF}(S^3)$.    }
\end{center}
\end{figure}

As another example, consider the right-handed trefoil $RHT$.  This knot has knot Floer homology given by 
\begin{displaymath}
\widehat{HFK}_j(RHT,i) = \left\{
	\begin{array}{lr}
		\mathbb{Z}_2\mathrm{\;if\;}(i,j) = (1,0)\\
		\mathbb{Z}_2\mathrm{\;if\;}(i,j) = (0,-1)\\
		\mathbb{Z}_2\mathrm{\;if\;}(i,j) = (-1,-2)\\
		0\mathrm{\;otherwise}
	\end{array}
	\right.
\end{displaymath}
Here, $i,j$ denote the Alexander and Maslov gradings, respectfully.  Since we know there is a spectral sequence, induced by the Alexander filtration on $\widehat{CF}(S^3)$, converging to $\widehat{HF}(S^3)\cong\mathbb{Z}$ (supported in grading 0), it follows that the $E^1$ page of this spectral sequence is given in Figure \ref{figure:tref}.

\begin{figure}
\labellist
\small\hair 2pt
\pinlabel $\mathbb{Z}_2$ at 47 43
\pinlabel $\mathbb{Z}_2$ at 47 13
\pinlabel $\mathbb{Z}_2$ at 15 43
\pinlabel $j$ at 16 3
\pinlabel $i$ at 2 13
\endlabellist
\begin{center}
\includegraphics{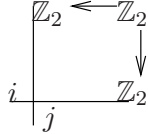}
\caption{  
\label{figure:ctref} A generating complex for the knot complex of the right-handed trefoil.   }
\end{center}
\end{figure}

By the symmetry of $CFK^\infty$, it follows that $CFK^\infty(RHT)$ is generated as a $\mathbb{Z}[U,U^{-1}]$--module by the complex in Figure \ref{figure:ctref}.

which shows that $d(S^3_1(RHT)) = -2$, a fact which more readily follows from Equation \ref{eq:alt}.

\begin{figure}
\labellist
\small\hair 2pt
\pinlabel $\mathbb{Z}^3_2$ at 33 35
\pinlabel $\mathbb{Z}_2$ at 61 66
\pinlabel $\mathbb{Z}_2$ at 5 11
\pinlabel Alexander at -33 37
\pinlabel Maslov at 35 78
\endlabellist
\begin{center}
\includegraphics{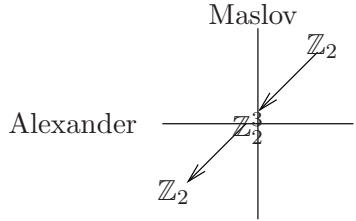}
\caption{\label{figure:4_1} The $E^1$ page of the spectral sequence $\widehat{HFK}(4_1)\Rightarrow\widehat{HF}(S^3)$. The markings on the arrows signify the ranks of the maps.   }
\end{center}
\end{figure}

Next, consider the figure eight knot, $4_1$.  This knot has knot Floer homology
\begin{displaymath}
\widehat{HFK}_j(4_1,i) = \left\{
	\begin{array}{lr}
		\mathbb{Z}_2\mathrm{\;if\;}(i,j) = (1,1)\\
		\mathbb{Z}_2^3\mathrm{\;if\;}(i,j) = (0,0)\\
		\mathbb{Z}_2\mathrm{\;if\;}(i,j) = (-1,-1)\\
		0\mathrm{\;otherwise}
	\end{array}\right.
\end{displaymath}
Again, by considering the spectral sequence $\widehat{HFK}(4_1)\Rightarrow\widehat{HF}(S^3)$, it follows that the $E^1$ page of this spectral sequence is given Figure \ref{figure:4_1}.

\begin{figure}
\labellist
\small\hair 2pt
\pinlabel $\mathbb{Z}_2$ at 52 52
\pinlabel $\mathbb{Z}_2$ at 52 24
\pinlabel $\mathbb{Z}_2$ at 24 24
\pinlabel $\mathbb{Z}_2$ at 24 52
\pinlabel $\mathbb{Z}_2$ at 14 14
\pinlabel $i$ at 3 24
\pinlabel $j$ at 23 5
\endlabellist
\begin{center}
\includegraphics{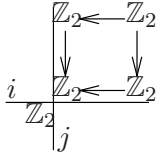}
\caption{\label{figure:c4_1} A generating complex for the figure eight knot, $4_1$.  Here we have two $\mathbb{Z}_2$'s at the origin---one is isolated while the other is part of a null-homologous ``box''.   }
\end{center}
\end{figure}

It follows that $CFK^\infty(4_1)$ is generated as a $\mathbb{Z}[U,U^{-1}]$--module by the complex shown in Figure \ref{figure:c4_1}.

\begin{figure}
\labellist
\small\hair 2pt
\pinlabel $\mathbb{Z}_2$ at 169 197
\pinlabel $\mathbb{Z}_2$ at 143 171
\pinlabel $\mathbb{Z}^3_2$ at 117 144
\pinlabel $\mathbb{Z}^2$ at 93 122
\pinlabel $\mathbb{Z}_2^3$ at 63 91
\pinlabel $\mathbb{Z}_2^3$ at 35 62
\pinlabel $\mathbb{Z}_2$ at 7 35
\pinlabel $1$ at 153 190
\pinlabel $2$ at 129 165
\pinlabel $1$ at 102 140
\pinlabel $1$ at 75 112
\pinlabel $2$ at 45 81
\pinlabel $1$ at 18 54
\pinlabel $\mathbb{Z}_2^3$ at 170 170
\pinlabel $\mathbb{Z}_2^3$ at 147 147
\pinlabel $\mathbb{Z}_2^3$ at 120 120
\pinlabel $\mathbb{Z}_2^3$ at 93 93
\pinlabel $\mathbb{Z}_2^3$ at 66 66
\pinlabel $\mathbb{Z}_2^3$ at 36 36
\pinlabel $\mathbb{Z}_2$ at 9 9
\pinlabel $1$ at 166 157
\pinlabel $2$ at 139 128
\pinlabel $1$ at 112 101
\pinlabel $2$ at 84 75
\pinlabel $2$ at 56 45
\pinlabel $1$ at 29 19
\endlabellist
\begin{center}
\includegraphics{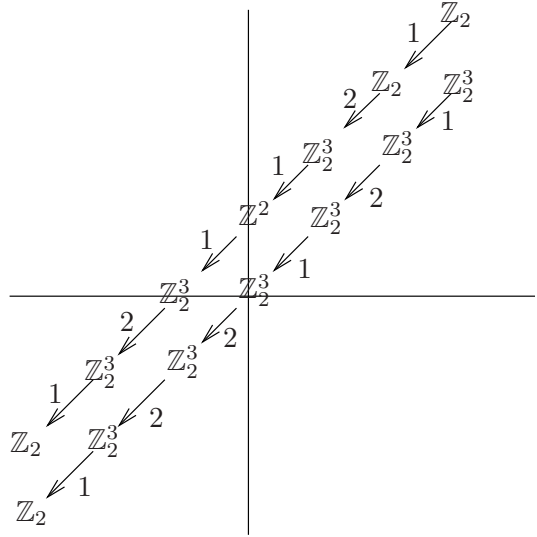}
\caption{\label{figure:C21} The $E^1$ term of the spectral sequence $\widehat{HFK}(C_{2,1}) \Rightarrow \widehat{HF}(S^3)$.   }
\end{center}
\end{figure}

Here, we see a single isolated $\mathbb{Z}_2$ at the origin plus a null-homologous ``box''.  It is then easy to see that $d(S^3_{\pm1}(4_1)) = 0$ (again, this follows more quickly from Equation \ref{eq:alt}).  

Recall that while the Kinoshita--Terasaka knot $K_{2,1}$ is smoothly slice, it is currently unknown if its Conway mutant $C_{2,1}$ is smoothly slice (though it is \emph{topologically} slice since it has trivial Alexander polynomial, by a result of Freedman \cite{Free,Free2}).  Indeed, we currently show that $d(S^3_1(C_{2,1})) = 0$, showing that our invariant gives no information.  In \cite{BG}, Baldwin and Gillam calculated the knot Floer homology polynomial\footnote{The knot Floer homology polynomial of a knot $K$ is defined to be $$\Sigma_{i,j}\dim_{\mathbb{Z}_2}\widehat{HFK}_j(S^3,K,i)\,q^jt^i$$.} of $C_{2,1}$ to be:
\begin{equation}\label{equation:c21}
(q^{-3}+q^{-2})t^{-3} + 3(q^{-2}+q^{-1})t^{-2} + 3(q^{-1}+1)t^{-1} + 3 + 2q + 3(q+q^2)t + 3(q^2+q^3)t^2 + (q^3+q^4)t^3
\end{equation}Similar to previous computations, it follows that the $E^1$ term of the spectral sequence  $\widehat{HFK}(C_{2,1}) \Rightarrow \widehat{HF}(S^3)$ is forced to be as in Figure \ref{figure:C21}.  From this, it follows that $CFK^\infty(C_{2,1})$  can be computed by a complex generated as a $\mathbb{Z}[U,U^{-1}]$--module with a single isolated $\mathbb{Z}_2$ at the origin plus a collection of null-homologous ``boxes''.  As in the computation for the figure eight, it follows that $d(S^3_{\pm1}(C_{2,1})) = 0$.


\subsection{An example session}
In this section we show an example session of our program \texttt{dCalc}.  We first input a generating complex for $CFK^\infty$ of the right-handed trefoil, as in Figure \ref{figure:ctref}.  We then form the complex for the connect-sum $RHT\#RHT$.  Finally we compute the correction terms of $S^3_{\pm1}(RHT\#RHT)$.

\begin{verbatim}
*****************************************************************
Welcome to the d invariant calculator!

This program computes the d invariants of +/-1 surgery on a knot

Copyright (C) 2009 Thomas Peters

This program is free software; you can redistribute it and/or
modify it under the terms of the GNU General Public License
as published by the Free Software Foundation; either version 2
of the license, or (at your option) any later version.

This program is distributed in the hope that it will be useful,
but WITHOUT ANY WARRANTY; without even the implied warranty of
MERCHANTABILITY or FITNESS FOR A PARTICULAR PURPOSE.  See the 
GNU General Public License for more details.

You should have received a copy of the GNU General Public License
along with this program; if not, write to the Free Software
Foundation, INC., 51 Franklin Street, Fifth Floor, Boston,
 MA 02110-1301, USA

email tpeters@math.columbia.edu with problems, bugs, etc
*****************************************************************
---------------------------------
   Main menu.  
(1) Enter a new knot
(2) View current knots
(3) Select a knot
(4) Connect-sum two knots
(0) Quit
----------------------------------
1
Enter the name of your knot
trefoil
Enter the knot vertex keys (non-neg integers).  input -1 to stop
0
1
2
-1
entered vertices 0,1,2,
enter the adjacency lists (type a vertex key, press enter,
	continue. input -1 to stop)
successors of 0:
1
2
-1
successors of 1:
-1
successors of 2:
-1
enter the bifiltration levels
(type i value, press enter, then type j value, then press enter)
(input -1 to stop)

F_i[0] = 1
F_j[0] = 1

F_i[1] = 0
F_j[1] = 1

F_i[2] = 1
F_j[2] = 0
added knot trefoil with adjacency list
[0]1,2,
[1]
[2]
and bifiltration levels
F(0) = (1,1)
F(1) = (0,1)
F(2) = (1,0)



---------------------------------
   Main menu.  
(1) Enter a new knot
(2) View current knots
(3) Select a knot
(4) Connect-sum two knots
(0) Quit
----------------------------------
4
Current knots are:
(index, name)
------------
(0, trefoil)
Enter the indices of the two knots to add
0
0
Computing tensor product...
Computation took 0min0sec.
Created knot trefoil#trefoil having adjacency list
[12]
[8]
[5]8,12,
[7]
[4]
[2]4,7,
[3]7,12,
[1]4,8,
[0]2,5,1,3,
and bifiltrations
F(12) = (2,0)
F(8) = (1,1)
F(5) = (2,1)
F(7) = (1,1)
F(4) = (0,2)
F(2) = (1,2)
F(3) = (2,1)
F(1) = (1,2)
F(0) = (2,2)
---------------------------------
   Main menu.  
(1) Enter a new knot
(2) View current knots
(3) Select a knot
(4) Connect-sum two knots
(0) Quit
----------------------------------
3
Current knots are:
(index, name)
------------
(0, trefoil)
(1, trefoil#trefoil)
input an index
1
What would you like to do with your knot complex?
(1) Print its adjacency list
(2) Show its bifiltration levels
(3) Check if it defines a complex
(4) Check if it is filtered
(5) Compute its homology
(6) Compute d invariants!
(7) Nothing--bring me back to the main menu
6
d(S^3_{+1}(K)) = -2
d(S^3_{-1}(K)) = 0
---------------------------------
   Main menu.  
(1) Enter a new knot
(2) View current knots
(3) Select a knot
(4) Connect-sum two knots
(0) Quit
----------------------------------
0
Really quit d calculator? (y/n) y
\end{verbatim}

\subsection{Issues with the implementation}
\label{section:issues}
\texttt{dCalc} does not do any checking on inputted complexes.  If one inputs a complex which does not come from a knot, \texttt{dCalc} may \texttt{return} \texttt{garbage} or have undefined behavior.  \texttt{dCalc} does, however, come with a few basic functions useful in determining the feasibility of a given complex.  For instance, it can check if the user's graph actually represents a complex. 

It is also worth mentioning that our implementation was for Floer homology with coefficients in $\mathbb{Z}_2$, so we are really computing correction terms for mod-2 coefficients.  It is an interesting question to determine whether or not $d$--invariants for Floer homology with $\mathbb{Z}_2$ coefficients can ever differ from $d$--invariants calculated with $\mathbb{Z}$ coefficients.

More seriously, by default \texttt{dCalc} uniquely identifies vertices by \texttt{int} keys.  One is therefore limited by the maximum value of \texttt{int}, \texttt{INT\_MAX} (this is defined in the header file \texttt{<limits.h>} and varies from platform to platform, though is guaranteed to be at least 32,767).  This is only realistically a problem after taking tensor products, where we rely on an explicit bijection $\mathbb{Z}\times\mathbb{Z}\to\mathbb{Z}$ to assign vertex keys for the tensor product.  This function is quadratic in its two arguments so it can grow quite quickly.  Surpassing \texttt{INT\_MAX} can result in undefined behavior (including segmentation faults).  If one were limited by this feature, one could change the underlying data structure of the vertex keys to a more flexible structure, for instance something like the \texttt{tuple} structure found in \texttt{python}, or to a larger integer structure, such as a \texttt{long} \texttt{unsigned} \texttt{int}.  The latter can be done by changing the line ``\texttt{typedef} \texttt{int} \texttt{KEYTYPE;}'' of \texttt{vertex.h} to, for instance, ``\texttt{typedef} \texttt{long} \texttt{KEYTYPE;}'' and them recompiling.  Of course, such operations increase run time.  One way to check if \texttt{INT\_MAX} has been exceeded (assuming \texttt{KEYTYPE} is not \texttt{unsigned}) is by printing out the adjacency matrix (or filtration levels) for a particular knot complex.  If negatives appear as keys, \texttt{INT\_MAX} has been surpassed (though, in principle, this need not be a necessary condition).

One place in which this program is inefficient memory-wise is in checking whether or not a given element in a complex is a boundary.  We do this by row-reduction.  If a complex has $n$ generators, the row-reduction requires a \texttt{char} array of roughly size $n^2$ to be allocated from the heap.  Of course, one should not need to create these matrices considering the homology itself can be checked just by performing an algorithm on the graph (see Baldwin and Gillam \cite{BG} for a discussion).  

We stress that in order to compute $d(S^3_{\pm1}(K))$ for a given knot, one must have at hand the filtered chain homotopy type of the $\mathbb{Z}[U,U^{-1}]$--module $CFK^{\infty}(K)$.  Computing these complexes is quite challenging, in general.  In the case that $K$ is alternating or is a torus knot, then one may recover $CFK^\infty(K)$ from the usually weaker invariant $\widehat{HFK}(K)$.  In the former case we have Equation \ref{eq:alt} and in the latter we have Ozsv\'ath--Szab\'o \cite[Theorem 1.2]{OSLens}, so we do not need to use a computer at all.  Depending on one's proficiency in Heegaard Floer homology, it is sometimes possible (though one should not expect in general) to calculate $CFK^{\infty}(K)$ from $\widehat{HFK}(K)$ (for instance, see the examples in Section \ref{section:examples}).

\bibliography{refs}
\bibliographystyle{amsalpha}

\end{document}